\documentclass[11pt]{amsart}

\usepackage{amssymb}
\usepackage{amscd}

\numberwithin{equation}{section}

\newtheorem{theorem}{Theorem}[section]

\newtheorem{lemma}[theorem]{Lemma}
\newtheorem{conjecture}[theorem]{Conjecture}
\newtheorem{remit}[theorem]{Remark}
\newtheorem{definit}[theorem]{Definition}

\newenvironment{remark}{\begin{remit}\rm}{\end{remit}}

\newcommand{\pp}{\mathbb{P}}

\newcommand{\zz}{\mathbb{Z}}

\newcommand{\cO}{\mathcal{O} }

\begin{document}

\title[Endomorphisms of hypersurfaces]{Endomorphisms of hypersurfaces of Fano manifolds of Picard number 1}
\date{}

\author{Insong Choe}
\address{Department of Mathematics, Konkuk University, 1 Hwayang-dong, Gwangjin-Gu, Seoul, Korea}
\email{ischoe@konkuk.ac.kr}

\thanks{}

\begin{abstract}
It is conjectured that a Fano manifold $X$ of Picard number 1 which
is not a projective space admits no endomorphisms of degree bigger
than 1. Beauville confirmed this for hypersurfaces of projective
space. We study this problem when $X$ is given by a hypersurface of
an arbitrary Fano manifold of Picard number 1.
\end{abstract}
\maketitle

In this paper, we are concerned with the following conjecture on the
endomorphisms of Fano manifolds of Picard number 1.
\begin{conjecture} \rm{(}\cite{AKP}, Conjecture 1.1\rm{)} Let $f$ be a surjective endomorphism of a Fano manifold $X$ of Picard number 1.
If $\deg f > 1$, then $X \cong \pp^n$.
\end{conjecture}
\noindent The reader is referred to the paper \cite{AKP} of Aprodu,
Kebekus, and Peternell for known results on this conjecture and
their recent enlargement of the list of Fano manifolds giving
evidence to this conjecture. Also one can see Fakhruddin's paper
\cite{F} for the arithmetic features of this problem.

On projective hypersurfaces, Beauville proved the following.
\begin{theorem} \rm{(}\cite{B}\rm{)}
A smooth complex projective hypersurface of dimension bigger than 1
and degree bigger than 2 admits no endomorphisms of degree bigger
than 1.
\end{theorem}
\noindent In particular, this shows that the projective
hypersurfaces of low degree, which are Fano manifolds of Picard
number 1, satisfy the conjectured property.

The proof of this result in \cite{B} relies mainly on the following
Chern number inequality which is obtained from the Hurwitz type
formula devised by Amerik, Rovinsky and Van de Ven.
\begin{lemma} \rm{(}\cite{ARV}, Corollary 1.2\rm{)} \label{ARV}
Let $f:X \to Y$ be a finite morphism between projective varieties of
dimension $n$. Let $L$ be a line bundle on $Y$ such that
$\Omega_Y(L)$ is globally generated. Then
\begin{equation} \label{ARVch}
\deg (f) \cdot c_n (\Omega_Y(L)) \le c_n (\Omega_X (f^*L)).
\end{equation}
\end{lemma}
\noindent To get the wanted result, Beauville combined this formula
with explicit computations of the Chern numbers of projective
hypersurfaces.

The aim of this paper is to apply the formula (\ref{ARVch}) in a
much more general context. The main idea is to replace the  explicit
Chern number computation in \cite{B} by another simple Chern number
inequality. Our result is the following.
\begin{theorem} \label{main}
Let $V$ be a Fano manifold with $Pic(V) \cong \zz$ generated by
$\cO_V(1)$. Assume that $\dim (V) \ge 4$ and $\Omega_V(l)$ is
globally generated. Let $X$ be a smooth hypersurface of $V$ cut out
by a member of  $|\cO_V(d)|$, $d \ge 2l$. Then $X$ admits no
surjective endomorphisms of degree bigger than 1.
\end{theorem}
An interesting special case is when $V$ is a {\it prime} Fano
manifold, which means by definition that the ample generator
$\cO_V(1)$ is very ample. For example, compact Hermitian symmetric
spaces and homogeneous contact manifolds of Picard number 1 are
known to be prime. In this case, $V$ is embedded in a projective
space $\pp^N$ via $ \cO_V(1)$ and thus $\Omega_V(2)$ is globally
generated, because it is a quotient of a globally generated bundle
$\Omega_{\pp^N} (2)$. Hence we have the following. Recall that the
{\it index} of $V$ is defined as the number $i$ satisfying $K_V
\cong \cO_V(-i)$.
\begin{theorem} \label{prime}
Let $V$ be a prime Fano manifold of Picard number 1 with index $i$.
Assume that $\dim V \ge 4$. Let $X$ be a smooth hypersurface of $V$
cut out by a member of $|\cO_V(d)|$, $4 \le d < i$. Then $X$ is a
Fano manifold of Picard number 1 and $X$ admits no surjective
endomorphisms of degree $>1$.
\end{theorem}
\noindent {\it Proof.} \ Non-existence of endomorphisms is a
consequence of the above theorem. So it suffices to check that $X$
is a Fano manifold of Picard number 1. From the Fano condition on
$V$, $H^1(V, \cO_V) = H^2(V, \cO_V)= 0$. By Lefschetz hyperplane
theorem, the restriction of divisors yields the isomorphism $Pic(X)
\cong Pic(V) \cong \zz$. The Fano condition on $X$ can be checked by
the adjunction formula
$$
K_X \cong K_V|_X \otimes [X]|_X \cong \cO_X(-i+d),
$$
where $\cO_X(1) = \cO_V(1)|_X$ is the ample generator of $Pic(X)$.
\qed

\vspace{0.1cm}
\section{Proof}
In this section, we prove Theorem \ref{main}. Suppose that $X$ and
$V$ satisfy the assumptions of Theorem \ref{main}. Let $\dim X = n
\ge 3$, let $h = c_1 (\cO_X(1))$. First we have an upper bound on a
Chern number.
\begin{lemma} \rm{(cf. \cite{B},  Proposition 1)} \  \label{upper}
If $X$ has a surjective endomorphism of degree bigger than 1, then
$$
c_n(\Omega_X(l)) \le l^n \cdot h^n.
$$
\end{lemma}
\noindent {\it Proof.} \ Let $f$ be an endomorphism of degree bigger
than 1. Let $m>1$ be the integer with $f^*\cO_X(1) \cong \cO_X(m)$.
Since $\Omega_V(l)$ is globally generated, $\Omega_X(l)$ is also
globally generated. Hence by Lemma \ref{ARV}, we have the following
inequality
\begin{equation} \label{upper2} m^n \cdot c_n(\Omega_X(l))
\le c_n(\Omega_X(lm)).
\end{equation}
Expanding $c_n(\Omega_X(lm))$ as a polynomial of $m$, we get
$$
c_n(\Omega_X(lm)) =  l^n \cdot h^n \cdot m^n + {(\rm  lower \ order
\ terms)}.
$$
By iterating the endomorphism $f$, we can make $m$ arbitrarily
large. Hence dividing both sides of (\ref{upper2}) by $m^n$ and
taking limit $m \to \infty$,
we get the wanted inequality.  \qed \\

Next from the conormal sequence twisted by $\cO_X(l)$:
$$
0 \to \cO_X(l-d) \to \Omega_V(l)|_X \to \Omega_X(l) \to 0,
$$
we get
$$
c(\Omega_X(l)) = c(\Omega_V(l)|_X) \cdot (1 + (l-d) h) ^{-1}.
$$
Hence the top Chern number has the following expression:
\begin{equation} \label{summation}
c_n(\Omega_X(l)) = \sum_{i=0}^{n} (d-l)^i \cdot
c_{n-i}(\Omega_V(l)|_X) \; h^i.
\end{equation}
Since $\Omega_V(l)$ is globally generated, its restriction to $X$ is
also globally generated. Hence the Chern classes
$c_i(\Omega_V(l)|_X)$ are represented by positive cycles of pure
codimension $i$ (cf. \cite{Fu}, Example 14.4.3). Hence for $d \ge
l$, each term appearing in the summation (\ref{summation}) is
nonnegative. In particular, we have
\begin{lemma}  \label{lower}
For $d \ge l$, $c_n(\Omega_X(l)) \ge (d-l)^n h^n$. \qed
\end{lemma}
To finish the proof of Theorem \ref{main}, we compare the bounds on
$c_n(\Omega_X(l))$ given by Lemma \ref{upper} and Lemma \ref{lower}
and get $d \le 2l$.

For $d=2l$, from Lemma \ref{upper} and (\ref{summation}),
$$
c_n(\Omega_X(l)) = \sum_{i=0}^{n} l^i \cdot c_{n-i}(\Omega_V(l)|_X)
\; h^i \ \le \ l^n h^n.
$$
Because $c_{n-i}(\Omega_V(l)|_X)  \cdot  h^i \ge 0$, we get
$c_n(\Omega_X(l)) = l^n h^n$ and
$$
c_{n-i}(\Omega_V(l)|_X)  \cdot  h^i= 0
$$
for each $i = 0, 1, \cdots, n-1$. But this is absurd: from the
vanishing $c_1(\Omega_V(l)|_X) = c_1(K_V \otimes \cO_V(l(n+1))|_X) =
0$, we get $K_V \cong \cO_V(-ln-l)$. From the well-known upper bound
on the index of Fano manifolds of Picard number 1 (\cite{KO} p.32),
we have $ln+l \le n+2$. Thus $l=1$ and $K_V \cong \cO_V(-n-1)$,
which implies that $V$ is a hyperquadric (\cite{Fuj} Theorem 2.3 or
\cite{BS} 3.1.6). But in this case it is known that $H^0(V,
\Omega_V(1)) = 0$ (\cite{S}, Section 4), contrary to the assumption
that $\Omega_V (l)$ is globally generated.

Therefore, we see that an existence of surjective endomorphism of
degree bigger than 1 yields $d < 2l$. \qed
\begin{remark} \label{final}
{\rm The above method doesn't work for the case $d \le l$. For
$l<d<2l$ (for $d=3$ when $V$ is prime), it is still possible that
the inequality in Lemma \ref{upper} yields a contradiction. For
example, the inequality is violated for cubic projective
hypersurfaces (\cite{B}). Also when $V$ is the Grassmannian $Gr(2,
4)$ and $X$ is a smooth hypersurface with $\cO_V(X) \cong \cO_V(3)$,
direct computation shows that $c_3(\Omega_X(2)) = 876$ which is
bigger than the expected upper bound $2^3 \cdot h^3$.}
\end{remark}
\begin{remark} \label{cohomology}
In the above proof, the Fano condition on $V$ was not essentially
used. In the proof, we needed the Fano condition only to guarantee
that $Pic(X) \cong \zz$ and it is generated by the hyperplane
section. Thus, the Fano condition on $V$ in Theorem \ref{main} can
be replaced by the following cohomological conditions:
$$
H^1(V, \cO_V) = H^2(V, \cO_V) = 0.
$$
\end{remark}

\section{Complete intersections}

In this section, we find a direct generalization of Theorem
\ref{main} to the case of complete intersections. Here we just
consider the case when $V$ is prime.
\begin{theorem}
Let $V$ be a prime Fano manifold of Picard number 1. Let $(d_i)_{1
\le i \le k}$ be positive numbers and let
$$
X = V \cap H_1 \cap H_2 \cap \cdots \cap H_k
$$
be a smooth subvariety of $V$ of codimension $k$ cut out by general
hypersurfaces $H_i \in |\cO_V(d_i)|$. Assume that $\dim X \ge 3$. If
$ d_i \ge 4$ for some $i$, then $X$ admits no surjective
endomorphisms of degree $>1$.
\end{theorem}
\noindent {\it Proof.} \ Since $H_i$ is general, we have a smooth
variety
$$
X_i := V \cap H_1 \cap \cdots \cap H_i
$$
for each $i$. By Lefschetz hyperplane theorem, $Pic(X_i) \cong \zz$
and it is generated by $\cO_V(1)|_{X_i}$. Rearranging the sequence
of hypersurfaces $\{ H_i \}$, we may assume that $d_k \ge 4$. The
smooth variety $X_{k-1}$ satisfies all the conditions of $V$ in
Theorem \ref{main}, possibly except for the Fano condition. Also, $X
= X_k$ is a smooth hypersurface of $X_{k-1}$ of degree $d_k \ge 4$.
Hence by Remark \ref{cohomology}, $X$ admits no surjective
endomorphisms. \qed \vspace{0.3cm}

{\it Acknowledgements.} \  We would like to thank J.-M. Hwang for
reading the first draft of this paper and giving a remark on the
prime condition in Theorem \ref{prime}. Also we would like to thank
the referee for informing the references \cite{BS}  and \cite{Fu}.

\end{document}